%
%
%
%
\input amstex.tex
\input amsppt.sty
\documentstyle{amsppt}
\advance\vsize2\baselineskip

\def\A{\Cal A}

\def\ro{\rho}
\def\bigint{\int}

\def\cal{\Cal}
\def\mathcal{\Cal}
\def\mathbb{\Bbb}

\def\U{\Cal U}

\def\R{\Bbb R}
\def\N{\Bbb N}

\def\phi{\varphi}
\def\e{\varepsilon}

\def\P{\Cal P}

\def\Q{\Bbb Q}

\def\exp{\operatorname{exp}}
\def\id{\operatorname{id}}

\def\Comp{\cal Comp}

\def\supp{\operatorname{supp}}

\def\min{\operatorname{min}}
\def\max{\operatorname{max}}
\def\inf{\operatorname{inf}}

\NoBlackBoxes \topmatter
\title  Hartman-Mycielski functor of non-metrizable compacta
\endtitle
\author Taras Radul and Du\v{s}an Repov\v{s}
\endauthor
\address
Department of Mechanics and Mathematics, Lviv National University,
 Universytetska st.,1, 79602 Lviv, Ukraine.\endaddress
 \email
tarasradul\@yahoo.co.uk
\endemail

\address Institute For Matematics, Physics and Mechanics, and 
University Of Ljubljana, P.O.Box 2964, Ljubljana, Slovenia 1001
\endaddress
\email dusan.repovs\@guest.arnes.si
\endemail
\thanks The research was supported by the Slovenian-Ukrainian grant
SLO-UKR 06-07/04.
\endthanks
\abstract We investigate some topological properties of a normal
functor $H$ introduced earlier by Radul which is a certain
functorial
compactification of the Hartman-Mycielski construction $HM$. We
show that $H$ is open and  find the condition when $HX$ is an
absolute retract  homeomorphic to the Tychonov cube.
\endabstract
\endtopmatter

\document

\centerline{\bf 1. Introduction}
\vskip 0.3cm

The general theory of functors acting on the category $\Comp$ of
compact Hausdorff spaces (compacta) and continuous mappings was
founded by Shchepin [Sh2]. He described some elementary
properties of such functors and defined the notion of the normal
functor which has become very fruitful. The classes  of all normal
and weakly normal functors include many classical constructions:
the hyperspace {\sl exp}, the space of probability measures $P$, the
superextension $\lambda$, the space of hyperspaces of inclusion
$G$, and many other functors (cf. [FZ] and [TZ]).

Let $X$ be a space and $d$  an  admissible metric on $X$ bounded
by $1$. By $HM(X)$ we shall denote the space of all  maps from
$[0,1)$ to the space~$X$ such that $f|[t_i,t_{i+1})\equiv{}${\sl const},
for some $0=t_0\le\dots\le t_n=1$, with respect to the  following metric $$
d_{HM}(f,g)=\bigint_0^1d(f(t),g(t))dt,\qquad f,g\in HM(X). $$

The construction of $HM(X)$ is known as the {\sl Hartman-Mycielski
construction} [HM].

For every $Z\in\Comp$ consider $$\align HM_n(Z)=\Bigl\{f\in
HM(Z)\mid &\text{ there exist } 0=t_1<\dots<t_{n+1}=1\\ &\text{
with }f|[t_i,t_{i+1})\equiv z_i\in Z, i=1,\dots,n\Bigr\}.
\endalign$$

Let $\U$ be the unique uniformity of $Z$. For every $U\in\U$ and
$\e>0$, let
$$<\alpha,U,\e>=\{\beta\in HM_n(Z)\mid m\{t\in
[0,1)\mid (\alpha(t),\beta(t'))\notin U\}<\e\}.$$
The sets
$<\alpha,U,\e>$ form a base of a compact Hausdorff topology in
$HM_n Z$.
Given a map $f:X\to Y$ in $\Comp$, define a map $HM_n X\to HM_n Y$
by the formula $HM_n F(\alpha)=f\circ\alpha$. Then $HM_n$ is a
normal functor in $\Comp$ (cf. [TZ; 2.5.2]).

For $X\in\Comp$ we consider the space $HMX$ with the topology
described above. In general, $HMX$ is not compact. Zarichnyi
 has asked if there exists a normal functor in $\Comp$ which contains
 all  functors $HM_n$ as subfunctors (see [TZ]).
Such a functor $H$ was constructed in [Ra]. It was shown in [RR]
that $HX$ is homeomorphic to the Hilbert cube for each
non-degenerated metrizable compactum $X$.

We investigate some topological properties of the space $HX$ for
non-metrizable compacta $X$. The main results of this paper are:

\proclaim {Theorem 1.1} $Hf$ is open if and only if $f$ is an open map.
\endproclaim

\proclaim {Theorem 1.2}  $HX$ is an absolute retract if and only if $X$ is
openly generated compactum of weight $\le\omega_1$.
\endproclaim

\proclaim {Theorem 1.3} $HX$ is homeomorphic to Tychonov cube if and only if
$X$ is  an openly generated $\chi$-homogeneous compactum of weight
$\omega_1$.
\endproclaim

\vskip 0.5cm

 \centerline{\bf 2. Construction of $H$ and its
connection with the functor of }
\centerline{\bf probability measures $P$}

\vskip 0.3cm

Let $X\in\Comp$. By $CX$ we denote the Banach space of all
continuous functions $\phi:X\to\R$ with the usual $\sup$-norm:
$\|\phi\| =\sup\{|\phi(x)|\mid x\in X\}$. We denote the segment
$[0,1]$ by $I$.

For $X\in\Comp$ let us define the uniformity of $HMX$. For each
$\phi\in C(X)$ and $a,b\in [0,1]$ with $a<b$ we define the
function $\phi_{(a,b)}:HMX\to \R$ by the following formula
$$\phi_{(a,b)}=\frac 1{(b-a)}\bigint_a^b\phi\circ\alpha(t)dt.$$
Define
$$S_{HM}(X)=\{\phi_{(a,b)}\mid \phi\in C(X) \ \hbox{and} \ (a,b)\subset [0,1)\}.$$

For $\phi_1,\dots,\phi_n\in S_{HM}(X)$ define a pseudometric
$\ro_{\phi_1,\dots,\phi_n}$ on $HMX$ by the formula
$$\ro_{\phi_1,\dots,\phi_n}(f,g)=\max\{|\phi_i(f)-\phi_i(g)|\mid i\in\{1,\dots,n\}\},$$
where $f,g\in HMX$.
The family of pseudometrics
$$\P=\{\ro_{\phi_1,\dots,\phi_n}\mid n\in\N, \ \hbox{where} \
\phi_1,\dots,\phi_n\in S_{HM}(X)\},$$ defines a totally bounded
uniformity $\U_{HMX}$ of $HMX$ (see [Ra]).

For each compactum $X$ we consider the uniform space
$(HX,\U_{HX})$ which is the completion of $(HMX,\U_{HMX})$ and the
topological space $HX$ with the topology induced by the uniformity
$\U_{HX}$. Since $\U_{HMX}$ is totally bounded, the space $HX$ is
compact.

Let $f:X\to Y$ be a continuous map. Define the map $HMf:HMX\to
HMY$ by the formula $HMf(\alpha)=f\circ\alpha$, for all $\alpha\in
HMX$. It was shown in [Ra] that the map
$HMf:(HMX,\U_{HMX})\to(HMY,\U_{HMY})$ is uniformly continuous.
Hence there exists the continuous map $Hf:HX\to HY$ such that
$Hf|HMX=HMf$. It is easy to see that $H:\Comp\to\Comp$ is a
covariant functor and $HM_n$ is a subfunctor of $H$ for each
$n\in\N$.

Let us remark that the family of functions $S_{HM}(X)$ embed $HMX$
in the product of closed intervals $\prod_{\phi_{(a,b)}\in
S_{HM}(X)}I_{\phi_{(a,b)}}$ where $I_{\phi_{(a,b)}}=[\min_{x\in X}
|\phi(x)|,\max_{x\in X}$ $|\phi(x)|]$. Thus, the space $HX$ is the
closure of the image of $HMX$.  We denote by
$p_{\phi_{(a,b)}}:HX\to I_{\phi_{(a,b)}}$ the restriction of the
natural projection.  Let us remark that the function $Hf$ could be
defined by the condition $p_{\phi_{(a,b)}}\circ Hf=p_{(\phi\circ
f)_{(a,b)}}$ for each $\phi_{(a,b)}\in S_{HM}(Y)$.

It is shown in [RR] that $HX$ is a convex subset of
$\prod_{\phi_{(a,b)}\in S_{HM}(X)}I_{\phi_{(a,b)}}$.

Define the map $e_1:HMX\times HMX\times I\to HMX$ by the condition
that $e_1(\alpha_1,\alpha_2,t)(l)$ is equal to $\alpha_1(l)$ if
$l<t$ and $\alpha_2(l)$ in the opposite case for
$\alpha_1,\alpha_2\in HMX$, $t\in I$ and $l\in [0,1)$. We consider
$HMX$ with the uniformity $\U_{HMX}$ and $I$ with the natural
metric. The map $e_1:HMX\times HMX\times I\to HMX$ is uniformly
continuous [RR].

Hence there exists the extension of $e_1$ to the continuous map
$e:HX\times HX\times I\to HX$. It is easy to check that
$e(\alpha,\alpha,t)=\alpha$, for each $\alpha\in HX$.

We recall that  $PX$ is the space of all nonnegative functionals
$\mu:C(X)\to\R$ with norm 1, and taken in the weak* topology for a
compactum $X$ (see [TZ] or [FZ] for more details). Recall that the
base of the weak* topology in $PX$ consists of the sets of the
form $O(\mu_0,f_1,\dots,f_n,\e)=\{\mu\in
PX\mid|\mu(f_i)-\mu_0(f_i)|<\e$ for every $1\le i\le n\}.$ Hence
we can consider $PX$ as a subspace of the product of closed
intervals $\prod_{\phi\in C(X)}I_\phi$ where $I_\phi=[\min_{x\in
X} |\phi(x)|,\max_{x\in X}$ $|\phi(x)|]$.  We denote by
$\pi_\phi:PX\to I_\phi$ the restriction of the natural projection.

For each $(a,b)\subset (0,1)$ we can define a map
$rX_{(a,b)}:HX\to PX$ by formula $\pi_\phi\circ
rX_{(a,b)}=p_{\phi_{(a,b)}}$. It is easy to check that
$rX_{(a,b)}$ is well defined, continuous and affine map.

We define as well a map $iX:PX\to HX$ by the formula
$p_{\phi_{(a,b)}}\circ iX=\pi_\phi$.  We have that
$rX_{(a,b)}\circ iX=\id_{PX}$, hence $rX_{(a,b)}$ is a retraction
for each $(a,b)\subset (0,1)$. The map $rX_{(0,1)}$ we denote
simply by $rX$.

Let us remark that $r_{(a,b)}:H\to P$ is a natural transformation.
(It means that for each map $f:X\to Y$ we have $Pf\circ
rX_{(a,b)}=rY_{(a,b)}\circ Hf$.) The same property is valid for
$i:P\to H$.

\vskip 0.5cm

 \centerline {\bf 3. Openess of the functor $H$}

\vskip 0.3cm

A subset $A\subset HX$ is called $e$-{\it convex} if
$e(\alpha,\beta,t)\in A$, for each $\alpha$, $\beta\in A$ and
$t\in I$. If, additionally, $A$ is convex, we call $A$ $H$-{\it
convex}.

We suppose that $f:X\to Y$ is  a continuous surjective map between
compacta during this section. The proofs of the next three lemmas
are easy checking on $HMX$ which is a dense subset of $HX$.

\proclaim {Lemma 3.1} For each $\mu$, $\nu\in HX$ and $t\in[0,1]$
we have $e(Hf(\mu),Hf(\nu),t)=Hf(e(\mu,\nu,t))$.
\endproclaim


\proclaim {Lemma 3.2} Consider any $\nu\in HX$ and $a$, $b$,
$c\in\R$ such that $0\le a<c<b\le 1$. Then we have
$p_{\phi_{(a,b)}}(\nu)=\frac{c-a}{b-a}p_{\phi_{(a,c)}}(\nu)+\frac{b-c}{b-a}p_{\phi_{(c,b)}}(\nu)$
for each $\nu\in HX$.
\endproclaim

\proclaim {Lemma 3.3} Let $t\in(0,1)$ and $(a,b)\subset (0,1)$.
For each $\mu$, $\nu\in HX$ and $\phi\in C(X)$ we have
$p_{\phi_{(a,b)}}(e(\mu,\nu,t))=p_{\phi_{(a,b)}}(\mu)$ if $b\le t$
and $p_{\phi_{(a,b)}}(e(\mu,\nu,t))=p_{\phi_{(a,b)}}(\nu)$ if
$t\le a$.
\endproclaim

\proclaim {Lemma 3.4} Let $A$ be a closed $H$-convex subset of
$HX$ and $\nu\notin A$. Then there exist $\phi\in C(X)$ and
$(a,b)\subset (0,1)$ such that.
$p_{\phi_{(a,b)}}(\nu)<p_{\phi_{(a,b)}}(\mu)$ for each $\mu\in A$.
\endproclaim

\demo {Proof} Suppose the contrary. We can for each $\mu\in A$
choose $\psi_\mu\in S_{HM}(X)$ such that
$p_{\psi_\mu}(\nu)<p_{\psi_\mu}(\mu)$. Since $A$ is compact, there
exist $\mu_1,\dots,\mu_n\in A$ such that for each $\mu\in A$ there
exists $i\in\{1,\dots,n\}$ such that
$p_{\psi_{\mu_i}}(\nu)<p_{\psi_{\mu_i}}(\mu)$. By Lemma 3.2 we can
choose a family of intervals $\{(a_i,b_i)\}_{i=1}^k$ such that
$b_i\le a_{i+1}$ and for each $i\in\{1,\dots,k\}$ a family of
function $\phi_{(a_i,b_i)}^1,\dots,\phi_{(a_i,b_i)}^{n_i}\in
S_{HM}(X)$ such that for each $\mu\in A$ there exist
$i\in\{1,\dots,k\}$ and $l\in\{1,\dots,n_i\}$ such that
$p_{\phi_{(a_i,b_i)}^l}(\nu)<p_{\phi_{(a_i,b_i)}^l}(\mu)$.

Consider the set $K=\{\mu\in A\mid p_{\phi_{(a_i,b_i)}^l}(\mu)\le
p_{\phi_{(a_i,b_i)}^l}(\nu)$ for each $i\in\{2,\dots,k\}$ and
$l\in\{1,\dots,n_i\}\}$. Then $K$ is a compact convex subset of
$A$ and for each $\mu\in A$ there exists $l\in\{1,\dots,n_1\}$
such that
$p_{\phi_{(a_1,b_1)}^l}(\nu)<p_{\phi_{(a_1,b_1)}^l}(\mu)$. Then
$rX_{(a_1,b_1)}(K)$ is a convex compact subset of $PX$ which
doesn't contain $rX_{(a_1,b_1)}(\nu)$. Then there exists
$\psi^1\in C(X)$ such that
$\pi_{\psi^1}(rX_{(a_1,b_1)}(\nu))<\pi_{\psi^1}(\eta)$ for each
$\eta\in rX_{(a_1,b_1)}(K)$. Hence, for each $\mu\in K$ we have
$p_{\psi^1_{(a_1,b_1)}}(\nu)<p_{\psi^1_{(a_1,b_1)}}(\mu)$.

Proceeding, we obtain $\psi^1,\dots,\psi^k\in C(X)$ such that for
each $\mu\in A$ there exists $i\in\{1,\dots,k\}$ such that
$p_{\psi^i_{(a_i,b_i)}}(\nu)<p_{\psi^i_{(a_i,b_i)}}(\mu)$.

By our supposition we can choose $\mu_i\in A$ for each
$i\in\{1,\dots,k\}$ such that $p_{\psi^i_{(a_i,b_i)}}(\mu_i)\le
p_{\psi^i_{(a_i,b_i)}}(\nu)$. Put $\xi_1=\mu_1$ and
$\xi_{i+1}=e(\xi_i,\mu_{i+1},b_i)$ for $i\in\{1,\dots,k-1\}$.
Since $A$ is $e$-convex, $\xi_k\in A$. By Lemma 3.3 we have
$p_{\psi^i_{(a_i,b_i)}}(\xi_k)\le p_{\psi^i_{(a_i,b_i)}}(\nu)$ for
each $i\in\{1,\dots,k\}$. Thus, we obtain a contradiction and the
lemma is proved.
\enddemo

The proof of the next lemma follows from Lemma 3.1 and the fact
that $Hf$ is affine map.

\proclaim {Lemma 3.5} $(Hf)^{-1}(\nu)$ is $H$-convex for each
$\nu\in HY$.
\endproclaim

Let $f:X\to Y$ be a map and $\phi\in C(X)$. By $\phi_*$ we denote
the function $\phi_*$:$Y\to \R$ defined by the formula
$\phi_*(y)=\inf(\phi(f^{-1}(y))$, $y\in Y$. It is known [DE] that
if $f$ is open then the function $\phi_*$ is continuous.

\demo {Proof of Theorem 1.1} Let $f:X\to Y$ be a map such that the
map $Hf:HX\to HY$ is open. Let us show that the map $Pf$ is open.
Consider any open set $U\subset PX$ and $\mu\in U$. Then
$(rX)^{-1}(U)$ is an open set in $HX$ and $iX(\mu)\in(rX)^{-1}(U)$
Since $Hf$ is an open map,$Hf((rX)^{-1}(U))$ is open in $HY$ and
$Hf(iX(\mu))\in Hf((rX)^{-1}(U))$. Since $r$ is a natural
transformation, we have $Hf((rX)^{-1}(U))\subset(rY)^{-1}(Pf(U))$.
We have $iY(Pf(\mu))=Hf(iX(\mu))$ or
$Pf(\mu)\in(iY)^{-1}(Hf((rX)^{-1}(U)))\subset(iY)^{-1}((rY)^{-1}(Pf(U)))=Pf(U)$.
 Since \newline $(iY)^{-1}(Hf((rX)^{-1}(U)))$ is open, the map $Pf$ is
 open. Hence $f$ is open as well [DE].

Now let a map $f:X\to Y$ be open. Let us suppose that $Hf$ is not
open. Then there exists $\mu_0\in HX$, a net $\{\nu_\alpha,
\alpha\in\A\}\subset O(Y)$ converging to $\nu_0=Hf(\mu_0)$ and a
neighborhood $W$ of $\mu_0$ such that $(Hf)^{-1}(\nu_\alpha)\cap
W=\emptyset$ for each $\alpha\in\A$. Since $HM(Y)$ is a dense
subset of $HY$, we can suppose that all $\nu_\alpha\in HM(Y)$ .
Since $HX$ is a compactum, we can assume that the net $A_\alpha=
(Hf)^{-1}(\nu_\alpha)$ converges in $\exp(HX)$ to some closed
subset $A\subset HX$. It is easy to check that $A\subset
(Hf)^{-1}(\nu_0)$ and $\mu_0\notin A$. By the Lemma 3.5 all the
sets $A_\alpha$ are $H$-convex. It is easy to see that $A$ is
$H$-convex as well. Since $\mu_0\notin A$, there exists by Lemma
3.4 $\phi\in C(X)$ and $(a,b)\subset (0,1)$ such that
$p_{\phi_{(a,b)}}(\mu_0)<p_{\phi_{(a,b)}}(\mu)$ for each $\mu\in
A$. Consider any $\alpha\in\A$. Let
$\{y_1,\dots,y_s\}=\nu_\alpha([0,1))$.
 Choose
for each $y_i$ the point $x_i$ such that $f(x_ i)= y_i$ and
$\phi(x_i)=\phi_*(y_i)$. Define a map
$j:\{y_1,\dots,y_s\}\to\{x_1,\dots,x_s\}$ by the formula
$j(y_i)=x_i$ and put $\mu_\alpha(t)=j\circ\nu_\alpha(t)$ for
$t\in[0,1)$. Let $\mu$ be a limit point of the net $\mu_\alpha$,
then $\mu\in A$. Since
$\phi_{(a,b)}(\mu_\alpha)=\phi_{*(a,b)}(\nu_\alpha)$, we have
$p_{\phi_{(a,b)}}(\mu)=p_{\phi_{*(a,b)}}(\nu_0)=p_{(\phi_*\circ
f)_{(a,b)}}(\mu_0)\le p_{\phi_{(a,b)}}(\mu_0)$. We have obtained
the contradiction and the theorem is proved.
\enddemo

\vfill\eject
 \centerline {\bf 4. Proofs}

\vskip 0.3cm

We will need some notations and facts from the theory of
non-metrizable compacta. See [10] for more details.

Let $\tau$ be an infinite cardinal number. A partially ordered set
$\A$ is called $\tau$-{\it complete}, if every subset of
cardinality $\le\tau$ has a least upper bound in $\A$. An inverse
system consisting of compacta and surjective bonding maps over a
$\tau$-complete indexing set is called $\tau$-complete. A
continuous $\tau$-complete system consisting of compacta of weight
$\le\tau$ is called a $\tau$-{\it system}.

As usual, by $\omega$ we denote the countable cardinal number.

A compactum $X$ is called {\it openly generated} if $X$ can be
represented as the limit of an $\omega$-system with open bonding
maps.

\demo {Proof of Theorem 1.2} It was shown in [RR] that $HX$ is
an absolute retract for each metrizable compactum $X$. So, we can
consider only non-metrizable case.

Let $X$ be an openly generated compactum of weight $\le\omega_1$.
By Theorem 1.1 the compactum $HX$ is also openly generated. Since
weight of $X$ (and $HX$ [Ra]) is  $\le\omega_1$, $HX$ is $AE(0)$.
Since $HX$ is a convex compactum, $HX$ is $AR$ [Fe].

Now, suppose $HX\in AR$. Since $rX:HX\to PX$ ia a retraction, $PX$
is an $AR$ too. Then $X$ is an openly generated compactum of
weight $\le\omega_1$[Fe]. The theorem is proved.
\enddemo

By $w(X)$ we denote the weight of the compactum $X$, by
$\chi(x,X)$ the character in the point $x$ and by $\chi(X)$ the
character of the space $X$. The space $X$ is called $\chi$-{\it
homogeneous} if for each $x,y\in X$ we have $\chi(x,X)=\chi(y,X)$.
We will use the following characterization of the Tychonov cube
$I^\tau$. An $AR$-compactum $X$ of weight $\tau$ is homeomorphic
to $I^\tau$ for an uncountable cardinal number $\tau$ if and only if $X$ is
$\chi$-homogeneous [Sh1].

Let $x\in X$. Define $\delta(x)\in HX$ by the condition
$p_{\phi_{(a,b)}}(\delta(x))=\phi(x)$ for each $\phi_{(a,b)}\in
S_{HM}(X)$.

\proclaim {Lemma 4.1} Let $f:X\to Y$ be an open map. Then $Hf$ has
a degenerate fiber if and only if $f$ has a degenerate fiber.
\endproclaim

\demo {Proof} Let $f:X\to Y$ be an open map such that there exists
$y\in Y$ with $f^{-1}(y)=\{x\}$, $x\in X$. Consider any $\mu\in
HX$ with $Hf(\mu)=\delta(y)$. Let us show that $\mu=\delta(x)$.
Consider any $\phi_{(a,b)}\in S_{HM}(X)$. Suppose that
$p_{\phi_{(a,b)}}(\mu)\ne\phi(x)$. We can assume that
$p_{\phi_{(a,b)}}(\mu)<\phi(x)$. By [Ra, Lemma 1] there exists a
function $\psi\in C(Y)$ such that $\psi(y)=\phi(x)$ and $\psi\circ
f\le\phi$. Then we have $p_{(\psi\circ f)_{(a,b)}}(\mu)\le
p_{\phi_{(a,b)}}(\mu)<\phi(x)$ and
$p_{\psi_{(a,b)}}(\delta(y))=p_{\psi_{(a,b)}}O(f)(\mu)=
p_{(\psi\circ f)_{(a,b)}}(\mu)<\phi(x)=\psi(y)$. Hence we obtain
the contradiction. Thus, $Hf$ has a degenerate fiber.

Now, suppose $f$ has no degenerate fiber. Consider any $\mu\in
HY$. Take any $y\in\supp \mu\subset Y$. Since $f$ is an open map
and $f^{-1}(y)$ is not a singleton, we can choose two closed
subsets $A_1,A_2\subset X$ such that $f(A_1)=f(A_2)=Y$ and
$(A_1\cap f^{-1}(y))\cap(A_2\cap f^{-1}(y)) =\emptyset$. Since the
functor $H$ preserves surjective maps [Ra], there exist $\mu_1\in
H(A_1)$ and $\mu_2\in H(A_2)$ such that $Hf(\mu_1)=Hf(\mu_2)=\mu$.
Since $y\in\supp\mu$, there exist $y_1\in\supp\mu_1\subset A_1$
and $y_2\in\supp\mu_2\subset A_2$ such that $f(y_1)=f(y_2)=y$.
Hence $\mu_1\ne\mu_2$ and the lemma is proved.
\enddemo

\proclaim {Lemma 4.2} An openly generated compactum $X$ of weight
$\omega_1$ is $\chi$-ho\-mo\-ge\-neo\-us if and only if $HX$ is
$\chi$-ho\-mo\-ge\-neous.
\endproclaim

\demo {Proof} Let $HX$ be $\chi$-homogeneous. Since the functor
$H$ preserves the weight [Ra], $HX$ is an absolute retract such
that $\chi(\mu,HX)=\omega_1$  for each $\mu\in HX$.  Then take any
$x\in X$ and suppose that there exists  $\{U_i\mid i\in\N\}$  a
countable base of open neighborhoods. Consider a family of
functions $\{\phi_i\in C(X)\mid i\in\N\}$ such that $\phi_i(x)=1$,
$\phi_i|X\setminus U_\i=0$. Then the family of function
$\{\phi_{i(a,b)}\mid i\in\N; a,b\in\Q\}$ define a countable base
of neighborhoods of $\delta(x)$ in $HX$. We obtain a
contradiction, hence $X$ is  $\chi$-homogeneous of $X$.

Now let $X$ be a $\chi$-homogeneous openly generated compactum of
weight $\omega_1$. Then $\chi(X) =\omega_1$ [Ra, Lemma 4]. Suppose
that there exists a point $\nu\in HX$ such that $\chi(\nu,
HX)<\omega_1$. Represent $X$ as the limit space of an
$\omega$-system $\{X_\alpha,p_\alpha,\A\}$ with open limit
projections $p_\alpha$. There exists $\alpha\in\A$ such that
$(Hp_\alpha)^{-1}(Hp_\alpha(\nu))=\{\nu\}$. By Lemma 4.2 there
exists a point $z\in X_\alpha$ such that $p_\alpha^{-1}(z)=\{x\}$,
$x\in X$. Hence $\chi(x,X)<\omega_1$ and we obtain the
contradiction. The lemma is proved.
\enddemo

\demo {Proof of Theorem 1.3}
 The proof of the theorem follows from Theorem 1.2 and Lemma 4.2.

\enddemo


\enddocument